\documentclass[a4paper,10pt]{amsart}

\usepackage{amssymb}
\usepackage{enumerate}
\usepackage{xcolor}
\usepackage{url}
\usepackage{hyperref}

\hypersetup{
colorlinks=true,
linkcolor=blue,
citecolor=blue,
urlcolor=blue
}

\newtheorem*{theorem*}{Theorem}
\newtheorem{lemma}{Lemma}
\newtheorem*{proposition*}{Proposition}

\theoremstyle{remark}
\newtheorem*{acknowledgements}{Acknowledgements}

\newcommand{\G}{\Gamma}
\newcommand{\dd}{\mathrm{d}}
\newcommand{\eps}{\varepsilon}

\begin{document}
\title{A sharp Mittag-Leffler bound for series with separated exponents}

\author[G.~Molteni]{Giuseppe Molteni}
\address[G.~Molteni]{Dipartimento di Matematica\\
         Universit\`{a} degli studi di Milano\\
         via Saldini 50\\
         I-20133 Milano\\
         Italy}
\email{giuseppe.molteni1@unimi.it}
\urladdr{\url{https://sites.unimi.it/molteni/}}

\keywords{Mittag-Leffler functions, extremal problems}
\subjclass[2020]{Primary 33E12, Secondary 26A33}


\date{\today.
}

\begin{abstract}
Let $\eta\geq 0$ and $\delta>0$, and let $(r_j)_{j\geq 0}$ be a sequence of real numbers satisfying $\eta
\leq r_0 < r_1 < r_2 < \cdots$ and $r_{j+1}-r_j \geq \delta$ for every $j\geq 0$. We prove that, for
$\delta\leq 1$ and every $x>0$,
\[
\sum_{j=0}^\infty \frac{x^{r_j}}{\G(1+r_j)}
\leq \sum_{j=0}^\infty \frac{x^{\eta+j\delta}}{\G(1+\eta+j\delta)}
=    x^\eta E_{\delta,1+\eta}(x^\delta),
\]
where $E_{\delta,1+\eta}$ denotes the two-parameter Mittag-Leffler function. Moreover, equality holds for
some $x>0$ if and only if $r_j = \eta + j\delta$ for every $j$. The proof combines a finite compression
argument for unimodal functions with a shift inequality for the Mittag-Leffler function. The upper bound
$\delta \leq 1$ is sharp, since counterexamples exist for every $1 < \delta < 2$.
\end{abstract}

\maketitle

\section{Introduction}\label{sec:1A}
The classical Mittag-Leffler function
\[
 E_\alpha(z)
 = \sum_{j=0}^{\infty}\frac{z^j}{\G(1+j\alpha)},
 \qquad
 \alpha > 0,
\]
and its two-parameter generalization
\[
 E_{\alpha,\beta}(z)
 = \sum_{j=0}^{\infty}\frac{z^j}{\G(\beta+j\alpha)},
 \qquad
 \alpha, \beta > 0,
\]
are natural variations of the exponential function and play an important role in fractional calculus and
the theory of special functions. Their analytic properties, including various monotonicity and
log-concavity properties, have been extensively investigated. For a general treatment,
see~\cite{GorenfloKilbasMainardiRogosin}. Examples of applications and related questions appear in \cite{
Abdelhakim,          
DeifOliveira,        
EdreiSaffVarga,      
GarrappaGerholdPopolizioSimon, 
TahaPecaricRoFarid}. 

In this paper, we consider a class of series extending the classical Mittag–Leffler functions, namely
series of the form
\begin{equation}\label{eq:0A}
 \sum_{j=0}^{\infty}\frac{x^{r_j}}{\G(1+r_j)},
\end{equation}
where the exponents $0\leq r_0<r_1<r_2<\cdots$ satisfy
\[
 r_{j+1}-r_j\geq\delta,
 \qquad
 j\geq 0,
\]
for some $\delta>0$. When $r_0=0$, the initial term is defined to be $1$ also at $x=0$, so that every
summand is continuous on $[0,+\infty)$ and $C^\infty$ on $(0,+\infty)$. The separation condition implies
that $r_j\geq r_0+j\delta$ for every $j$, so that the uniform convergence on every compact subset of
$[0,+\infty)$ immediately follows from Stirling's formula. This proves the continuity on $[0,+\infty)$ of
the resulting function, and the analogous argument for the series of derivatives of every order proves
that the function is $C^\infty$ on $(0,+\infty)$. Finally, only finitely many exponents $r_j$ are lower
than $1$, hence the function is also absolutely continuous on every interval $[0,R]$.

For fixed $x>0$, $\eta\geq 0$ and $\delta>0$, we seek to maximize~\eqref{eq:0A} over all sequences
satisfying $r_0\geq \eta$ and the above separation condition. Our main result identifies the arithmetic
progression $r_j=\eta + j\delta$ as the unique maximizing sequence when $0<\delta\leq 1$. More precisely,
we prove the following theorem.
\begin{theorem*}
Let $\eta\geq 0$ and $0<\delta\leq 1$. If $0\leq \eta \leq r_0 < r_1 < r_2 < \cdots$ satisfies $r_{j+1} -
r_j \geq\delta$ for every $j\geq 0$, then for every $x>0$,
\[
 \sum_{j=0}^{\infty}\frac{x^{r_j}}{\Gamma(1+r_j)}
 \leq \sum_{j=0}^{\infty}\frac{x^{\eta+j\delta}}{\Gamma(1+\eta+j\delta)}
 =    x^\eta E_{\delta,1+\eta}(x^\delta).
\]
Equality holds for some $x>0$ if and only if $r_j=\eta+j\delta$ for every $j\geq 0$ (in which case
equality holds for every $x>0$).
\end{theorem*}
Thus, a `gentle' modification of a Mittag-Leffler function provides a sharp upper bound for a whole
class of series with arbitrarily spaced exponents. In the special case $\eta=0$ and $\delta=1$, the
theorem produces a corresponding extremal property of the exponential function:
\[
 \sum_{j=0}^{\infty}\frac{x^{r_j}}{\G(1+r_j)}
 \leq e^x
\]
for every $x>0$, with equality for some positive $x$ precisely when the exponents are the nonnegative
integers, in which case equality holds for every $x>0$.

The proof combines two ingredients. The first is a finite compression principle for unimodal functions:
among finite configurations of points on the nonnegative real axis with a prescribed minimum separation,
the sum of a unimodal function can be increased, or at least preserved, by reducing the gaps until all
points form an arithmetic progression. The second ingredient is the fact that, for $0<\delta \leq 1$ and
every fixed $x>0$, the function
\[
 t\mapsto \sum_{j=0}^{\infty}\frac{x^{t+j\delta}}{\G(1+t+j\delta)}
 = x^t E_{\delta,1+t}(x^\delta)
\]
is strictly decreasing on $[0,+\infty)$.

The shift inequality is related to the complete monotonicity properties of Mittag-Leffler functions
studied by Simon~\cite{Simon}, as well as to the log-concavity results of Gerhold and
Simon~\cite{GerholdSimon}. The present argument uses these analytic properties to obtain an extremal
result for arbitrary separated sequences of exponents, rather than only for the standard Mittag-Leffler
series.

Finally, the restriction $\delta\leq 1$ is sharp: for every $1<\delta<2$ and every $0\leq\eta<\delta-1$,
there exists an admissible sequence for which the asserted inequality is reversed for all sufficiently
large $x$. Indeed, choose $\rho\in(\eta,\delta-1)$ and set $r_j=\rho+j\delta$ for every $j\geq 0$. We
show that
\[
 \sum_{j=0}^{\infty}\frac{x^{r_j}}{\Gamma(1 + r_j)}
 > \sum_{j=0}^{\infty}\frac{x^{\eta + j\delta}}{\Gamma(1 + \eta + j\delta)}
\]
for all sufficiently large $x$.

\begin{acknowledgements}
The author is a member of the INdAM group GNSAGA.
\end{acknowledgements}

\newpage
\section{Proof}\label{sec:2A}
We begin with an elementary lemma.
\begin{lemma}\label{lem:1A} %
Let $f\colon [0,\infty)\to[0,\infty)$ be unimodal. Thus, there exists $c\geq 0$ such that $f$ is
nondecreasing on $[0,c]$ and nonincreasing on $[c,\infty)$. Let $\eta\geq 0$, $\delta>0$ and let $0\leq
\eta \leq x_0 < x_1 < \cdots <x_N$ satisfy $x_{j+1} - x_j \geq \delta$ for every $j$. Then there exists
$t$ with $\eta \leq t\leq x_N - N\delta$ such that
\[
 \sum_{j=0}^N f(x_j)
 \leq \sum_{j=0}^N f(t+\delta j).
\]
\end{lemma}
\begin{proof}
By replacing $f(x)$ with $g(x):= f(\eta+\delta x)$ it suffices to consider the case $\eta=0$ and $\delta=1$.\\
Suppose that, for some $k\in\{0,\ldots,N-1\}$, $x_{k+1} - x_k = 1+\eps$ for some $\eps>0$. We show that
this gap can be reduced to $1$ without decreasing the value of $\sum_{j=0}^N f(x_j)$, while all the other
gaps remain unchanged.\\
Recall that $c$ is a mode of $f$, so that $f$ is nondecreasing on $[0,c]$ and nonincreasing on
$[c,\infty)$. There are three cases.
\medskip\\
\textbf{Case 1: $x_k\geq c$.} %
Then $x_j \geq c$ for every $j \geq k+1$. Move the right-hand block to the left by $\eps$, setting
\[
y_j
:= \begin{cases}
   x_j      & 0  \leq j\leq k,  \\
   x_j-\eps & k+1\leq j\leq N.
  \end{cases}
\]
Then $y_{k+1} - y_k = x_{k+1} - \eps - x_k = 1$ and all the other gaps are unchanged. Moreover, for
$j\geq k+1$, $y_j = x_j - \eps \geq x_{k+1}- \eps = x_k + 1 \geq c$. Since $y_j\leq x_j$ and $f$ is
nonincreasing on $[c,\infty)$, $f(y_j) \geq f(x_j)$. Hence
\[
 \sum_{j=0}^N f(y_j)\geq \sum_{j=0}^N f(x_j).
\]
\medskip\\
\textbf{Case 2: $x_{k+1}\leq c$.} %
Move the left-hand block to the right by $\eps$, setting
\[
y_j
:= \begin{cases}
   x_j+\eps & 0  \leq j\leq k,\\
   x_j      & k+1\leq j\leq N.
  \end{cases}
\]
Again, $y_{k+1} - y_k = x_{k+1} - x_k - \eps =1$, and all other gaps remain unchanged. Furthermore, for
$j\leq k$, $y_j\leq x_k+\eps = x_{k+1} - 1 < c$. Since $y_j\geq x_j$ and $f$ is nondecreasing on $[0,c]$,
$f(y_j) \geq f(x_j)$. Thus the sum does not decrease.
\medskip\\
\textbf{Case 3: $x_k<c<x_{k+1}$.} %
Set $ a:=c - x_k$ and $b:=x_{k+1} - c$. Then $a,b>0$ and $a + b = 1 + \eps$. Since $\eps < a + b$, we can
choose $\alpha,\beta\geq 0$ such that $\alpha + \beta = \eps$ and $\alpha \leq a$, $\beta \leq b$.
%
Now move the left-hand block to the right by $\alpha$ and the right-hand block to the left by $\beta$,
setting
\[
y_j
:= \begin{cases}
   x_j+\alpha & 0  \leq j\leq k,  \\
   x_j-\beta  & k+1\leq j\leq N.
  \end{cases}
\]
The central gap becomes $y_{k+1} - y_k = (x_{k+1}-\beta) - (x_k+\alpha)
= 1$.
The other gaps are unchanged.

For $j\leq k$, we have
\[
 x_j
 \leq y_j
 =    x_j+\alpha
 \leq x_k+\alpha
 \leq x_k+a
 =    c.
\]
Hence both $x_j$ and $y_j$ belong to $[0,c]$, and, since $f$ is nondecreasing on this interval,
\[
 f(y_j)\geq f(x_j).
\]
Similarly, for $j\geq k+1$,
\[
 c=   x_{k+1}-b
 \leq x_{k+1}-\beta
 \leq x_j-\beta
 =    y_j
 \leq x_j.
\]
Thus both $y_j$ and $x_j$ belong to $[c,\infty)$, and the monotonicity of $f$ there gives
\[
 f(y_j)\geq f(x_j).
\]
Thus again the total sum does not decrease.\\
We have proved that whenever one of the gaps is strictly larger than $1$, it can be reduced to $1$,
without changing the other gaps. Repeating this operation, after finitely many steps all gaps are equal
to $1$. Hence the resulting configuration has the form
\[
 t,
 t+1,
 \ldots,
 t+N
\]
for some $t\geq 0$, and
\[
 \sum_{j=0}^N f(x_j)
 \leq
 \sum_{j=0}^N f(t+j).
\]
Rescaling by $\delta$ and shifting by $\eta$ proves the general version of the lemma. Finally, we notice
that at every step the rightmost point never increases, so the final rightmost point satisfies
$t+N\delta\leq x_N$. This proves the bound for $t$.
\end{proof}
We shall also need a strict version.
\begin{lemma}\label{lem:2A}
Under the hypotheses of Lemma~\ref{lem:1A}, suppose additionally that $f$ is strictly increasing on
$[0,c]$ and strictly decreasing on $[c,\infty)$. If at least one gap exceeds $\delta$, then the
compression procedure in the proof of Lemma~\ref{lem:1A} produces $t\in [\eta,x_N-N\delta]$ such that
\[
 \sum_{j=0}^N f(x_j)
 < \sum_{j=0}^N f(t+j\delta).
\]
\end{lemma}
\begin{proof}
If one of the gaps is strictly larger than $\delta$, the corresponding compression step in the proof of
Lemma~\ref{lem:1A} moves at least one point strictly toward the unique maximizer $c$. Strict monotonicity
therefore makes the corresponding value of $f$ strictly larger. Hence at least one compression step is
strict, and equality is impossible.
\end{proof}
We next establish the property needed to apply the previous lemmas.
\begin{lemma}\label{lem:3A}
For every $x>0$, the function
\[
 f_x(t) := \frac{x^t}{\G(1+t)},
 \quad
 t\geq 0,
\]
is strictly unimodal. More precisely, there exists $c\geq 0$ such that $f_x$ is strictly increasing on
$[0,c]$ and strictly decreasing on $[c,\infty)$. Let $\gamma$ be the Euler--Mascheroni constant. When
$x\leq \exp(-\gamma)=0.56\ldots$ the mode is $c=0$ and $f_x$ is strictly decreasing on $[0,+\infty)$.
\end{lemma}
\begin{proof}
The function $f_x$ is positive in the given range and its logarithmic derivative is
\[
 (\ln f_x)'(t) = \ln x - \psi(1+t),
\]
where $\psi$ denotes the digamma function $\psi:=\Gamma'/\Gamma$. It is strictly increasing, satisfies
$\psi(1)=-\gamma$, and tends to $+\infty$ (see~\cite[Theorem~1.2.5] {AnAsRoy}). Thus $(\ln f_x)'$ has at
most one zero. Moreover, $f_x$ is eventually decreasing and $f_x(t)\to 0$ as $t\to\infty$ by Stirling's
formula. Hence, if a zero of $(\ln f_x)'$ exists, it is the unique point where $f_x$ attains its maximum.
\end{proof}
For $\delta>0$, $t\geq 0$ and $x\geq 0$, define
\[
 S_{t,\delta}(x)
 := \sum_{j=0}^\infty \frac{x^{t+j\delta}}{\G(1+t+j\delta)}
  = x^t E_{\delta,1+t}(x^\delta).
\]
In particular,
\[
 S_{0,\delta}(x)
 = \sum_{j=0}^{\infty} \frac{x^{j\delta}}{\G(1+j\delta)}
 = E_\delta(x^\delta).
\]
The next lemma is the key analytic input.
\begin{lemma}\label{lem:4A}
Let $0<\delta\leq 1$. For every fixed $x>0$, the function $t\mapsto S_{t,\delta}(x)$ is strictly
decreasing on $[0,\infty)$.
\end{lemma}
\begin{proof}
Fix $\delta$. We first prove that
\begin{equation}\label{eq:1A}
 S_{0,\delta}'(x)\ge S_{0,\delta}(x),
 \qquad x > 0,
\end{equation}
where the derivative is with respect to $x$. %
When $\delta=1$, this holds as equality since $S_{0,1}(x) = e^x$. Assume $0<\delta<1$. A representation
due to Simon~\cite[Eq.~(3.4)]{Simon} gives
\begin{equation}\label{eq:2A}
 S_{0,\delta}(x)
 = \frac{e^x}{\delta} - \varphi_\delta(x),
\end{equation}
where
\[
 \varphi_\delta(x)
 = \frac{\sin(\pi\delta)}{\pi} \int_0^\infty \frac{v^{\delta-1}e^{-xv}}{v^{2\delta}-2\cos(\pi\delta)v^\delta+1}\dd v.
\]
Thus, $\varphi_\delta$ is the Laplace transform of a positive measure, hence it is completely monotone
(as remarked in~\cite[Theorem~1.1]{Simon}). Differentiating~\eqref{eq:2A} and subtracting~\eqref{eq:2A}
yields
\[
 S'_{0,\delta}(x) - S_{0,\delta}(x)
 = \varphi_\delta(x) - \varphi'_\delta(x)
\]
and~\eqref{eq:1A} follows, since $\varphi_\delta(x)\geq 0$ and $\varphi'_\delta(x)\leq 0$.\\
We next prove that, for every $t\geq 0$,
\begin{equation}\label{eq:3A}
 S_{t,\delta}'(x)\ge S_{t,\delta}(x),
 \qquad x > 0,
 \quad  t\geq 0,
\end{equation}
where, again, the derivative is with respect to $x$. %
For this purpose we exploit an integral formula connecting $S_{t,\delta}$ to $S_{0,\delta}$. Assume
$t>0$. Then the beta integral (see Definition~1.1.3 and Theorem~1.1.4 in \cite{AnAsRoy}) gives the
equality
\[
 \int_0^x(x-u)^{t-1}u^{j\delta}\dd u
 = x^{t+j\delta}B(t,1+j\delta)
 = x^{t+j\delta}\frac{\G(t)\G(1+j\delta)}{\G(1+t+j\delta)},
\]
so that
\begin{align}
 S_{t,\delta}(x)
 &= \sum_{j=0}^{\infty}\frac{x^{t+j\delta}}{\G(1+t+j\delta)}
  = \frac{1}{\G(t)}\sum_{j=0}^{\infty}\frac{1}{\G(1+j\delta)}\int_0^x(x-u)^{t-1}u^{j\delta}\dd u \notag\\
 &= \frac{1}{\G(t)}\int_0^x (x-u)^{t-1}S_{0,\delta}(u)\dd u,                                     \label{eq:4A}
\end{align}
where the interchange of summation and integration is justified by Tonelli's theorem.
An integration by parts of~\eqref{eq:4A} yields
\[
 S_{t,\delta}(x)
 = \frac{x^t}{\Gamma(1+t)}
  + \frac{1}{\Gamma(1+t)} \int_0^x (x-u)^t S_{0,\delta}'(u)\dd u.
\]
Differentiating this identity and subtracting~\eqref{eq:4A} we obtain
\[
 S_{t,\delta}'(x)-S_{t,\delta}(x)
 =  \frac{x^{t-1}}{\Gamma(t)}
  + \frac{1}{\Gamma(t)} \int_0^x (x-u)^{t-1} \bigl(S_{0,\delta}'(u)-S_{0,\delta}(u)\bigr)\dd u.
\]
This formula and~\eqref{eq:1A} prove~\eqref{eq:3A} for $t>0$. When $t=0$, \eqref{eq:3A} coincides
with~\eqref{eq:1A}, so \eqref{eq:3A} is fully proved.\\
It follows from~\eqref{eq:3A} that $x\mapsto e^{-x}S_{t,\delta}(x)$ is nondecreasing for every $t\geq 0$,
and hence
\begin{equation}\label{eq:5A}
 S_{t,\delta}(u)\leq e^{-(x-u)}S_{t,\delta}(x),
 \qquad 0\leq u\leq x.
\end{equation}
Let $h>0$ and $t\geq 0$. The equality in~\eqref{eq:4A} can be generalized to
\begin{equation}\label{eq:6A}
 S_{t+h,\delta}(x)
 = \frac{1}{\G(h)}\int_0^x (x-u)^{h-1}S_{t,\delta}(u)\dd u
\end{equation}
by the same beta-integral argument used to obtain~\eqref{eq:4A}, or by applying the semigroup property of
the Riemann–Liouville fractional integral to~\eqref{eq:4A} (see~\cite[Chapter~2, Theorem~2.2]{Diethelm}).
Note that for $t=0$, it is simply~\eqref{eq:4A} with $t$ replaced by $h$ (so, \eqref{eq:6A} truly
generalizes~\eqref{eq:4A}).
Combining~\eqref{eq:6A} and~\eqref{eq:5A}, we obtain
\begin{align*}
 S_{t+h,\delta}(x)
 &\leq \frac{S_{t,\delta}(x)}{\Gamma(h)}\int_0^x (x-u)^{h-1}e^{-(x-u)}\dd u\\
 &=    \frac{S_{t,\delta}(x)}{\Gamma(h)}\int_0^x v^{h-1}e^{-v}\dd v
 <     S_{t,\delta}(x),
\end{align*}
where the last step comes from the inequality $\int_0^x v^{h-1}e^{-v}\dd v < \int_0^\infty v^{h-1}e^{-v}
\dd v = \G(h)$ and the positivity of $S_{t,\delta}(x)$. This proves the strict monotonicity on
$[0,\infty)$.
\end{proof}
\begin{proof}[Proof of the theorem]
Fix $x> 0$ and $N\geq 0$. By Lemma \ref{lem:3A}, the function $f_x(t) = x^t/\G(1+t)$ is strictly
unimodal. Since $r_{j+1} - r_j\geq \delta$ for every $j$, Lemma~\ref{lem:1A} gives a number $t_N$ in
$[\eta,r_N - N\delta]$ such that
\[
 \sum_{j=0}^N f_x(r_j)
 \leq \sum_{j=0}^N f_x(t_N+j\delta)
 \leq S_{t_N,\delta}(x).
\]
Lemma~\ref{lem:4A} gives $S_{t_N,\delta}(x) < S_{\eta,\delta}(x)$ in case $t_N > \eta$, and evidently
$S_{t_N,\delta}(x)=S_{\eta,\delta}(x)$ if $t_N=\eta$. In either case,
\[
 \sum_{j=0}^N f_x(r_j)
 \leq S_{\eta,\delta}(x).
\]
Letting $N\to\infty$ yields
\[
 \sum_{j=0}^\infty \frac{x^{r_j}}{\G(1+r_j)}
 =    \sum_{j=0}^\infty f_x(r_j)
 \leq S_{\eta,\delta}(x)
 =    \sum_{j=0}^\infty \frac{x^{\eta+j\delta}}{\G(1+\eta+j\delta)},
\]
which is the first part of the theorem. It remains to prove that equality is possible only when $r_j =
\eta + j\delta$ for every $j$.\\
Suppose some gap is strictly larger than $\delta$, and choose $N$ so that this gap lies among $r_0,
\ldots, r_N$. Apply the compression algorithm proving Lemma~\ref{lem:1A} to this block. Its strict form
in Lemma~\ref{lem:2A} gives $t_N\geq \eta$ with
\[
 D:= \sum_{j=0}^N\big(f_x(t_N+j\delta)-f_x(r_j)\big)
  >  0.
\]
From Lemma~\ref{lem:1A} we also know that $t_N+N\delta\leq r_N$. Define
\[
 r'_j
 := \begin{cases}
        t_N+j\delta  & \text{if $0\leq j\leq N$},\\
        r_j          & \text{if $      j >   N$}.
    \end{cases}
\]
The modified sequence has $r'_0 = t_N\geq \eta$ and all its gaps are at least $\delta$, including the
boundary gap, since $r'_{N+1} - r'_N = r_{N+1} - (t_N + N\delta) \geq r_{N+1} - r_N \geq \delta$. Thus,
the already established non-strict theorem applies to the sequence $(r'_j)$ and shows that
$\sum_{j=0}^\infty f_x(r'_j) \leq S_{\eta,\delta}(x)$. On the other hand, the tails of sequences $(r_j)$
and $(r'_j)$ coincide, so that a subtraction gives the equality
\[
 \sum_{j=0}^\infty f_x(r'_j)
 = \sum_{j=0}^\infty f_x(r_j) + D.
\]
Hence
\[
 \sum_{j=0}^\infty f_x(r_j)
 =    \sum_{j=0}^\infty f_x(r'_j) - D
 \leq S_{\eta,\delta}(x) - D
 <    S_{\eta,\delta}(x),
\]
and the strict inequality is proved. If no gap exceeds $\delta$, then $r_j = r_0 + j\delta$ for every $j$
and Lemma~\ref{lem:4A} gives strictness when $r_0>\eta$; when $r_0=\eta$, equality follows from the
definition of $S_{\eta,\delta}$.
\end{proof}
The following proposition shows that the monotonicity stated in Lemma~\ref{lem:4A} may fail when $\delta
> 1$, at least when $x$ is large enough.
\begin{proposition*}
Let $1 < \delta < 2$. Choose $\eta \in [0,\delta-1)$ and $\rho \in (\eta,\delta-1)$. Then
\[
 S_{\rho,\delta}(x) > S_{\eta,\delta}(x)
\]
for all sufficiently large $x$.
\end{proposition*}
\begin{proof}
We use the standard asymptotic expansion of the two-parameter Mittag-Leffler function; see, for example,
\cite[Theorem~4.3 and Proposition~4.5]{GorenfloKilbasMainardiRogosin}. Since $S_{\rho,\delta}(x) = x^\rho
E_{\delta,1+\rho}(x^\delta)$, we deduce that for $1 < \delta < 2$ and $0 \leq \eta < \rho < \delta-1$:
\begin{align*}
S_{\rho,\delta}(x)
 &= \frac{e^x}{\delta}
    - \frac{x^{\rho-\delta}}{\Gamma(1+\rho-\delta)}
    + o(x^{\rho-\delta}),\\
S_{\eta,\delta}(x)
 &= \frac{e^x}{\delta} + O(x^{\eta-\delta})
  = \frac{e^x}{\delta} + o(x^{\rho-\delta}),
\end{align*}
as $x\to\infty$. Thus, the exponentially large terms in $S_{\rho,\delta}(x)$ and $S_{\eta,\delta}(x)$
are equal, and
\[
 S_{\rho,\delta}(x) - S_{\eta,\delta}(x)
 = -\frac{x^{\rho-\delta}}{\G(1 + \rho - \delta)} + o(x^{\rho-\delta})
\]
as $x$ tends to $\infty$. The constant $\G(1 + \rho - \delta)$ is negative, because $-1 < 1+\rho-\delta < 0$.
Therefore
\[
 S_{\rho,\delta}(x) - S_{\eta,\delta}(x) > 0
\]
for all sufficiently large $x$.
\end{proof}
\noindent %
This gives an explicit counterexample to the extension of the theorem beyond $\delta=1$: if $1<\delta<2$,
pick any $\eta\in [0,\delta-1)$, $\rho\in(\eta,\delta-1)$ and take $r_j = \rho + j\delta$ for every $j$.
Then $r_0=\rho > \eta$ and $r_{j+1} - r_j = \delta$, but, for sufficiently large $x$,
\[
 \sum_{j=0}^\infty \frac{x^{r_j}}{\G(1+r_j)}
 =  S_{\rho,\delta}(x)
 >  S_{\eta,\delta}(x)
 =  \sum_{j=0}^\infty \frac{x^{\eta + j\delta}}{\G(1+\eta+j\delta)}.
\]


\end{document}